\documentclass{amsart}

\usepackage{amsmath}
\usepackage{amssymb}
\usepackage{graphicx}
\usepackage{booktabs}
\usepackage{makecell}

\renewcommand{\S}{\mathbb S}
\renewcommand{\H}{\mathbb H}
\newcommand{\CP}{\mathbb C\mathrm P}
\newcommand{\CH}{\mathbb C\mathrm H}
\newcommand{\Z}{\mathbb Z}
\newcommand{\R}{\mathbb R}

\newcommand{\set}[1]{ \left\{ #1 \right\} }
\newcommand{\SO}{\mathrm{SO}}
\renewcommand{\O}{\mathrm{O}}
\newcommand{\U}{\mathrm U}

\newcommand{\SL}{\mathrm{SL}}

\newtheorem{theorem}{Theorem}[section]

\theoremstyle{definition}

\theoremstyle{remark}
\newtheorem{remark}[theorem]{Remark}

\numberwithin{equation}{section}

\title[4-dimensional Thurston geometries]{Thurston geometries in dimension four\\from a Riemannian perspective}

\author{Marie D'haene}
\address{KU Leuven, Department of Mathematics, Celestijnenlaan 200B--Box 2400, 3001 Leuven, Belgium}
\email{marie.dhaene@kuleuven.be}
\thanks{The author is supported by Methusalem grant METH/21/03–long term structural funding of the Flemish Government.}

\begin{document}

\begin{abstract}
	In this survey we focus on a special class of homogeneous manifolds called Thur\-ston geometries.
	We give special attention to the four-dimensional Thurston geometries with 4 or 5-dimensional isometry group which are not a product (except for $\mathrm F^4$). 
	These are the manifolds $\mathrm{Sol}^4_0$, $\mathrm{Sol}^4_1$, $\mathrm{Sol}^4_{m,n}$ and $\mathrm{Nil}^4$.
	We give a description of each of the spaces and we exhibit all Riemannian metrics which are invariant under the action of their associated group.
\end{abstract}

\maketitle

Riemannian homogeneous spaces are a fundamental class of manifolds whose study requires techniques from geometry, algebra and group theory.
Thurston geometries are a subset of Riemannian homogeneous spaces, famous due to Thurston’s geometrization for 3-manifolds where these appear as fundamental building blocks.
This geometrization theorem concerns compact orientable 3-manifolds and it is very similar to the uniformization theorem for compact orientable surfaces.
Namely, it says that any 3-manifold can be cut up into pieces along 2-spheres and 2-tori, and the interior of each of these pieces admits a Riemannian metric which is locally isometric to one of eight 3-dimensional Thurston geometries.
These are $\R^3$, $\S^3$, $\H^3$, $\S^2 \times \R$,  $\H^2\times\R$, $\widetilde{\SL}(2,\R)$, $\mathrm{Nil}^3$ and $\mathrm{Sol}^3$ (see Table~\ref{tab:3-geometries}).
\begin{table}[t]
	\centering
	\begin{tabular}{c c c}
		\toprule
		Thurston geometry & Stabilizer & dimension of $G$ \\
		\midrule
		$\R^3,\; \S^3,\; \H^3$
		& $\SO(3)$ & 6 \\[1ex]
		$\S^2\times\R,\; \H^2\times\R, \; \widetilde{\SL}(2,\R), \; \mathrm{Nil}^3$
		& $\SO(2)$ & 5 \\[1ex]
		$\mathrm{Sol}^3$
		& $\set{1}$ & 3 \\
		\bottomrule \\
	\end{tabular}
	\caption{
	Based on~\cite{maier}.	
	Three-dimensional Thurston geometries ordered by their stabilizer of the \emph{connected component of the identity} of the group~$G$ that acts on them.
	}
	\label{tab:3-geometries}
\end{table}
Hence understanding these model spaces is an essential part of the study of 3-manifolds and as such they have already received a lot of attention in the literature.
In the context of isometrically immersed submanifolds, we refer to the following selection of works: Daniel~\cite{daniel2007isometric} studied the fundamental equations for surfaces in 3-dimensional homogeneous spaces, which include $\S^2\times\R$, $\H^2\times\R$, $\widetilde{\SL}(2,\R)$ and $\mathrm{Nil}^3$. 
Parallel and totally geodesic surfaces in these spaces were classified in~\cite{Belkhelfa-Dillen-Inoguchi}, while the classification in all remaining three-dimensional homogeneous spaces, including~$\text{Sol}^3$, can be found in \cite{I-VdV-2007} and \cite{I-VdV-2008}.
Moreover, totally geodesic and parallel submanifolds of the space forms were covered in~\cite{Lawson}, \cite{Takeuchi} and \cite{Backes-Reckziegel}.

In this article, however, the main focus is on the 4-dimensional Thurston geometries, which were classified in the PhD thesis of Filipkiewicz~\cite{Filip} (see Table~\ref{tab:4-geometries}).
\begin{table}
	\centering
	\begin{tabular}{c c c}
		\toprule
		Thurston geometry & Stabilizer & dimension of $G$ \\
		\midrule
		$\R^4,\; \S^4,\; \H^4$
		& $\SO(4)$ & 10 \\[1ex]
		$\CP^2,\; \CH^2$
		& $\U(2)$ & 9 \\[1ex]
		$\S^3 \times \R,\; \H^3 \times \R$
		& $\SO(3)$ & 7 \\[1ex]
		\makecell{$\S^2 \times \R^2,\; \S^2 \times \S^2$ \\ $\H^2 \times \S^2,\; \H^2 \times \R^2,\; \H^2 \times \H^2$}
		& $\SO(2) \times \SO(2)$ & 6 \\[2.5ex]
		$\text{Sol}^4_0,\; \mathrm F^4, \; \widetilde{\SL}(2,\R) \times \R,\; \text{Nil}^3 \times \R$
		& $\SO(2)$ & 5 \\[1ex]
		$\text{Nil}^4,\; \text{Sol}^4_1,\; \text{Sol}^4_{m,n}$
		& $\set{1}$ & 4 \\
		\bottomrule \\
	\end{tabular}
	\caption{
	Based on~\cite{Wall}.	
	Four-dimensional Thurston geometries ordered by their stabilizer of the \emph{connected component of the identity} of the group~$G$ that acts on them.
	In $\mathrm{Sol}^4_{m,n}$, $m$ and $n$ are positive integers and when $m=n$, $\mathrm{Sol}^4_{m,n} \cong \mathrm{Sol}^3 \times \R$.
	}
	\label{tab:4-geometries}
\end{table}
Although no geometrization is known in dimension~4, the study of 4-dimensional geometries and their submanifolds is interesting on its own from a Riemannian geometric point of view.
It provides essential insight in the structure of these spaces.
In many of the 4-dimensional geometries, more specifically the real and complex space forms and products of lower-dimensional geometries, ample research has been done regarding their isometrically immersed submanifolds.
Hence in this article we focus on other 4-dimensional geometries, $\mathrm{Sol}^4_0$, $\mathrm{Sol}^4_1$, $\mathrm{Sol}^4_{m,n}$ and $\mathrm{Nil}^4$.
Note that there is a remaining Thurston geometry, called $\mathrm{F}^4$ which does not appear in our present study, as this space requires different techniques than the ones mentioned earlier and hence does not fit in the scope of this article. 

An important remark is that the references mentioned above study Thurston geometries but, more often than not, they are not approached as such.
By this we mean that there is no systematic study preceding the study of submanifolds in these spaces regarding the family of Riemannian metrics which is invariant under the action of the corresponding group. 
The aim of the present article is to provide such a framework for the aforementioned subset of 4-dimensional Thurston spaces, mentioned in Table~\ref{tab:geometries-groups} with the associated groups.
\begin{table}
	\centering
	\begin{tabular}{c c}
		\toprule
		Thurston geometry & $G$ \\
		\midrule
		$\mathrm{Sol}^4_0$ & $\mathrm{Sol}^4_0 \rtimes (\O(2)\times \Z/2\Z)$ \\[1ex]
		$\mathrm{Sol}^4_1$ & $\mathrm{Sol}^4_1 \rtimes D_4$ \\[1ex] 
		$\mathrm{Sol}^4_{m,n} (m\ne n)$ & $\mathrm{Sol}^4_{m,n} \rtimes (\Z/2\Z)^3$ \\[1ex]
		$\mathrm{Nil}^4$ & $\mathrm{Nil}^4 \rtimes (\Z/2\Z)^2$ \\[1ex]
		$\mathrm F^4$ & $\R^2 \rtimes \SL^\pm(2,\R)$ \\ 
		\bottomrule \\
	\end{tabular}
	\caption{
	Based on \cite{hillman2002four} and \cite{wall-proceedings}.	
	A selection of 4-dimensional Thurston geometries with 4- or 5-dimensional isometry group.
	In the second column we provide the group $G$ that act on the manifolds in the first column.
	The manifolds together with the group form Thurston geometries.
	The group $\SL^\pm(2,\R)$ consists of matrices of determinant $+1$ or $-1$.
	}
	\label{tab:geometries-groups}
\end{table}
The results we present are not deep, but they are essential for the study of submanifolds.

\section{Thurston geometries}

The origin of Thurston geometries can be traced back to a very fundamental question: ``How does the topology of a surface relate to its geometry and vice versa?''
In this context, geometry refers to a Riemannian structure.

An example of the interplay between geometry and topology is the famous theorem of Gauss and Bonnet.
Namely, for a surface $S$ we have
$\int_S K \, dA = 2\pi \chi(S)$, relating the total Gaussian curvature on the left-hand side to a multiple of the Euler characteristic $\chi(S)$, a purely topological quantity, on the right-hand side.
This result shows that the topology of a surface prescribes its curvature and vice versa.
Another theorem that links the topology of a surface to its geometry is the theorem of Bonnet and Myers, which states that a Riemannian manifold with sectional curvature (a generalization of the Gaussian curvature of a surface) bounded from below by a positive constant must be compact and have finite first fundamental group, which again highlights the connection between curvature and topology.

The theorems mentioned above are very strong results, however, the most beautiful interplay between geometry and topology is illustrated by the theorem of Klein and Poincaré: any compact orientable surface admits a Riemannian metric with constant sectional curvature $k$.
Moreover, if we denote by $g$ the genus of $S$, then $k > 0$ if and only if $g = 0$, $k = 0$ if and only if
$g = 1$ and $k < 0$ if and only if $g > 2$.
This is often referred to as the uniformization theorem
of compact orientable surfaces.
A nice corollary of this statement is that any such surface
admits a Riemannian universal covering by $\R^2$ if $g = 1$, by $\S^2$ if $g = 0$ or by $\H^2$ if $g > 2$.
In this sense we say that the Riemannian manifolds $\R^2$, $\S^2$ and $\H^2$ are models for surfaces.
Upon encountering this theorem, one is immediately inclined to wonder if a similar result holds for 3-manifolds.
However, the geometry of 3-manifolds is much richer and more complicated than that of surfaces, since we easily see that not every compact orientable 3-manifold can be endowed with a metric of constant sectional curvature, for example, one can easily show that $\S^2\times\S^1$ cannot be endowed with a constant curvature metric.
Therefore, the list of 3-dimensional model geometries includes more spaces than $\R^3$, $\S^3$ and $\H^3$.
In an attempt to state an analogue of the 2-dimensional geometrization for 3-manifolds, Thurston pinpointed the special properties of the spaces $\R^2$, $\S^2$ and $\H^2$ and the definition of a Thurston geometry was born.
Following \cite{thurston-book} and \cite{bonahon}, a Thurston geometry is a pair $(X, G)$ of a smooth manifold $X$ and a Lie group $G \subset \mathrm{Diff}(X)$ satisfying the following properties.
\begin{itemize}
	\item[(a)] $X$ is connected and simply connected.
	\item[(b)] $G$ acts transitively on $X$ with compact stabilizers. Note that this implies the existence of a $G$-invariant Riemannian metric on $X$.
	\item[(c)] $G$ is maximal, in the sense that there does not exist a larger Lie group containing $G$ that also acts transitively on $X$ with compact stabilizers. 
	\item[(d)] There exists a finite volume manifold which is modelled by the pair $(X,G)$, meaning that there exists a finite volume Riemannian manifold $(M,g)$ that admits $(X,\tilde g)$ as a Riemannian universal cover (here $\tilde g$ is any $G$-invariant Riemannian metric on $X$).
\end{itemize}
\begin{remark}
	The definition of Thurston in~\cite{thurston-book} requires the existence of a \emph{compact} manifold modelled by the pair $(X,G)$, while Filipkiewicz~\cite{Filip} relaxes this to \emph{finite volume}.
	For 3-dimensional geometries, changing the \emph{finite volume} condition to compactness does not alter the spaces satisfying the definition of a Thurston geometry, as all 3-dimensional Thurston geometries are models for compact manifolds.
	However, in the 4-dimensional case, making this change does alter the list, in the sense that the space $\mathrm{F}^4$ only models finite volume manifolds.
	This is the only 4-dimensional geometry that does not model any compact manifold.
\end{remark}
There are exactly eight pairs $(X,G)$ with $\dim X = 3$ that satisfy the definition above, as Thurston shows in~\cite{thurston-book}.
These comprise the eight 3-dimensional Thurston geometries, listed in Table~\ref{tab:3-geometries}.
They are an essential part of the 3-dimensional analogue of the uniformization of surfaces, which is called the Thurston geometrization conjecture.
The proof is in large part due to Thurston and was finished in 2003 by Perelman, making it the first Millennium Prize Problem to be solved.
Indeed, the Poincaré conjecture is a consequence of the Thurston geometrization conjecture and given that the Poincaré conjecture had been an open problem for a hundred years, the proof of Perelman was an enormous breakthrough in geometry.
Furthermore, the Thurston geometrization conjecture gives a deep understanding of the interplay between the topology and geometry of 3-manifolds.

Four-manifolds are not as well understood as 3-manifolds, which is reflected in the lack of geometrization-type theorem for 4-manifolds.
However, in an attempt to understand these, we still study 4-dimensional Thurston geometries.
As mentioned earlier, they are pinpointed by Filipkiewicz~\cite{Filip}, and can be found in Table~\ref{tab:4-geometries}.

\section{A Riemannian description of some Thurston geometries}
Since Thurston geometries are Riemannian homogeneous spaces enjoying many symmetries, they are an interesting class of manifolds in which to study isometrically immersed submanifolds.
The first step in this process is to determine all Riemannian metrics which are compatible with the group structure, i.e.\ for a pair $(X,G)$ we have to determine the family of $G$-invariant metrics.
Only after we have accomplished this, we can study submanifolds of a Thurston geometry.
Up until today, no systematic study of $G$-invariant metrics has been performed in the non-product Thurston 4-spaces with 4- or 5-dimensional isometry group (see the last two rows of Table~\ref{tab:4-geometries}).
These are the spaces $\mathrm{Sol}^4_0$, $\mathrm F^4$, $\mathrm{Sol}^4_{m,n}$, $\mathrm{Sol}^4_1$ and $\mathrm{Nil}^4$, which can be found in Table~\ref{tab:geometries-groups} along with the corresponding group.
In this section we treat each geometry separately (except $\mathrm{F}^4$), describing its geometry.
We refer to~\cite{hillman2002four},~\cite{wall-proceedings} and \cite{Filip} for all spaces except $\mathrm{Sol}^4_0$, for which we refer to previous work of the author~\cite{sol40}. 

\subsection*{$\boldsymbol{\mathrm{Sol}^4_0}$}
The manifold associated to this Thurston geometry is a solvable Lie group containing the following matrices:
\[
 \begin{pmatrix}
	e^t & 0 & 0 & x \\ 
	0 & e^t & 0 & y \\
	0 & 0 & e^{-2t} & z \\
	0 & 0 & 0 & 1
 \end{pmatrix}
\]
where $t,x,y,z \in \R$.
The Lie algebra is spanned by the basis $\set{e_1,e_2,e_3,e_4}$ satisfying the following commutation relations:
\begin{equation*}
    \begin{aligned}
        &[e_1,e_2] = e_2, \quad [e_1,e_3] = e_3, \quad [e_1,e_4] = -2e_4, \\
        &[e_2,e_3] = [e_2,e_4] = [e_3,e_4] = 0.
    \end{aligned}
\end{equation*}
The left invariant vector fields determined by
$e_1$, $e_2$, $e_3$ and $e_4$ are
\begin{equation*}
	E_1=\frac{\partial}{\partial t},
	\quad
	E_2=e^{t}\frac{\partial}{\partial x},
	\quad
	E_3=e^{t}\frac{\partial}{\partial y},
	\quad
	E_4=e^{-2t}\frac{\partial}{\partial z}.
\end{equation*}
The group associated to this Thurston geometry is $G=\mathrm{Sol}^4_0 \rtimes (\O(2) \times \Z/2\Z)$ where $\mathrm{Sol}^4_0$ acts by left translation, $\O(2)$ acts in the $xy$-plane and $\Z/2\Z$ acts by reflecting the $z$-coordinate.

To determine all $G$-invariant Riemannian metrics on the manifold, we note that $\mathrm{Sol}^4_0 = G/K$ is a \emph{reductive} homogeneous space, where $K$ denotes the stabilizer. 
This means that there exists a subspace $\mathfrak p$ of the Lie algebra $\mathfrak g$ of $G$ such that $\mathfrak g = \mathfrak k \oplus \mathfrak p$ with the property that $\mathfrak p$ is invariant under the adjoint action of $K$.
For a reductive homogeneous space there is a bijection between $\mathrm{Ad}_K$-invariant inner products on~$\mathfrak p$ and $G$-invariant metrics on $G/K$.
Using this correspondence we find that, up to homothety, there is a unique $G$-invariant Riemannian metric on $\mathrm{Sol}^4_0$ given by 
\begin{equation*}
	dt^2 + e^{-2t}(dx^2+dy^2) + e^{4t}dz^2.
\end{equation*}
Note that the frame $\set{E_1,E_2,E_3,E_4}$ is orthonormal with respect to this metric.

As mentioned by Wall~\cite{wall-proceedings}, $\mathrm{Sol}^4_0$ admits a compatible complex structure.
In fact, it is unique up to isomorphism of the group $G$ and the structure is moreover integrable, making $\mathrm{Sol}^4_0$ a complex manifold (see~\cite{sol40}).
Note however that it is not Kähler, but one can show that the metric given above is globally conformal to a Kähler metric (the conformal factor is $e^{2t}$).

\subsection*{$\boldsymbol{\mathrm{Sol}^4_{m,n}}$}
The manifold of this Thurston geometry is also a solvable Lie group consisting of the matrices 
\[
\begin{pmatrix}
	e^{at} & 0 & 0 & x \\
	0 & e^{bt} & 0 & y \\
	0 & 0 & e^{ct} & z \\
	0 & 0 & 0 & 1
\end{pmatrix}	
\]
where $t,x,y,z \in \R$ and $a < b < c$ such that $e^{at}$, $e^{bt}$ and $e^{ct}$ are the roots of the polynomial $x^3 - m x^2 + nx - 1$.
Note that $m$ and $n$ are positive integers such that the roots are all \emph{real} and \emph{distinct}.
The associated Lie algebra is spanned by $\set{e_1,e_2,e_3,e_4}$ satisfying the commutation relations 
\begin{equation*}
    \begin{aligned}
        &[e_1,e_2] = a e_2, \quad [e_1,e_3] = b e_3, \quad [e_1,e_4] = c e_4, \\
        &[e_2,e_3] = [e_2,e_4] = [e_3,e_4] = 0.
    \end{aligned}
\end{equation*}
The left invariant vector fields associated to $e_1$, $e_2$, $e_3$ and $e_4$ are
\begin{equation*}
	E_1= \frac{\partial}{\partial t},
	\quad
	E_2=e^{at} \frac{\partial}{\partial x},
	\quad
	E_3=e^{bt} \frac{\partial}{\partial y},
	\quad
	E_4=e^{ct}\frac{\partial}{\partial z}.
\end{equation*}
Note that if $m = n$, then one of the roots equals 1 and $\mathrm{Sol}^4_{m,n}$ can be associated with $\mathrm{Sol}^3 \times \R$. 
We will not consider this case here, since we are only concerned with non-product spaces.
Hence from now on we will assume that $m$ and $n$ are different.
In that case, the associated group is $\mathrm{Sol}^4_{m,n} \rtimes (\Z/2\Z)^3$. 
Each $\Z/2\Z$ factor is generated by reflecting either the $x$-, $y$- or $z$-coordinate.
Another important remark to make is that if we allow $m$ and $n$ to be such that two of the roots are equal, we can identify $\mathrm{Sol}^4_{m,n}$ with $\mathrm{Sol}^4_0$.
In this case there is an additional symmetry in the $x$- and $y$-coordinates, corresponding to the $\O(2)$ factor in the group associated to $\mathrm{Sol}^4_0$.

Also $\mathrm{Sol}^4_{m,n}$ is a \emph{reductive} homogeneous space and using the same technique as for $\mathrm{Sol}^4_0$, we find that with $m \ne n$ there is a unique metric up to homothety which is invariant under the action of the associated group, given by
\begin{equation*}
	dt^2 + e^{-2at}dx^2 + e^{-2bt} dy^2 + e^{-2ct}dz^2.
\end{equation*}
Note that the frame $\set{E_1,E_2,E_3,E_4}$ is orthonormal with respect to this metric.

In~\cite{wall-proceedings} it is shown that $\mathrm{Sol}^4_{m,n}$ does not admit any compatible complex structures.

\subsection*{$\boldsymbol{\mathrm{Sol}^4_1}$}
The underlying manifold of this Thurston geometry is again a solvable Lie group and it consists of matrices of the form
\[
  M(t,x,y,z) =
	\begin{pmatrix}
        1 & x & z\\
        0 & t & y\\
        0 & 0 & 1
    \end{pmatrix}
\]
where $t,x,y,z\in \R$ and $t>0$.
The Lie algebra is spanned by the basis $\set{e_1,e_2,e_3,e_4}$ satisfying the following commutation relations:
\begin{equation*}
    \begin{aligned}
        &[e_1,e_2] = -e_2, \quad [e_1,e_3] = e_3, \quad [e_2,e_3] = e_4, \\
        &[e_1,e_4] = [e_2,e_4] = [e_3,e_4] = 0.
    \end{aligned}
\end{equation*}
The left invariant vector fields determined by
$e_1$, $e_2$, $e_3$ and $e_4$ are
\begin{equation*}
	E_1=t \frac{\partial}{\partial t}+x \frac{\partial}{\partial x},
	\quad
	E_2=\frac{\partial}{\partial x},
	\quad
	E_3=t\frac{\partial}{\partial y}+x \frac{\partial}{\partial z},
	\quad
	E_4=\frac{\partial}{\partial z}.
\end{equation*}
The group that acts on $\mathrm{Sol}^4_1$ making it a Thurston geometry is $\mathrm{Sol}^4_1 \rtimes D_4$ (where $D_4$ denotes the dihedral group of 8 elements).
$\mathrm{Sol}^4_1$ acts by left translations while the action of $D_4$ is given by the action of its generators:
\begin{equation*}
    s: M(t,x,y,z)\mapsto M(t,x,-y,-z), \quad   
    r: M(t,x,y,z) \mapsto M(1/t,-y/t,x/t,z-xy/t).
\end{equation*}
Again using that $\mathrm{Sol}^4_1$ is a reductive homogeneous space, we find, up to homothety, a two-parameter family of $G$-invariant metrics given by
\begin{equation*}
	\frac{1}{t^2} \Big( 
		(x^2+\tau_1) dt^2 + t^2 dx^2 + (1 + \tau_2 x^2) dy^2 + \tau_2 t^2 dz^2 - 2tx (dt\,dx + \tau_2 dy\,dz)
	\Big)
\end{equation*}
where $\tau_1,\tau_2\in \R^+$.
Note that the frame $\set{\frac{E_1}{\sqrt{\tau_1}},E_2,E_3,\frac{E_4}{\sqrt{\tau_2}}}$ is orthonormal with respect to this metric.

In~\cite{wall-proceedings} it is shown that $\mathrm{Sol}^4_1$ admits a compatible complex structure, and there are two isomorphism classes.
It is not Kähler. 

\subsection*{$\boldsymbol{\mathrm{Nil}^4}$}
The manifold associated to this Thurston geometry is a nilpotent Lie group and is most easily described as a semi-direct product $\R \ltimes_\theta \R^3$ where for each $t \in \R$ we have the following automorphism of $\R^3$:
\[
\theta(t) =
\begin{pmatrix}
	1 & t & \frac{t^2}{2} \\
	0 & 1 & t \\
	0 & 0 & 1
\end{pmatrix} \in \mathrm{Nil}^3.
\]
The Lie algebra is spanned by the basis $\set{e_1,e_2,e_3,e_4}$ with the following brackets:
\begin{equation*}
    \begin{aligned}
        &[e_1,e_2] = 0, \quad [e_1,e_3] = e_2, \quad [e_2,e_3] = e_3, \\
        &[e_2,e_3] = [e_2,e_4] = [e_3,e_4] = 0.
    \end{aligned}
\end{equation*}
The left invariant vector fields determined by
$e_1$, $e_2$, $e_3$ and $e_4$ are
\begin{equation*}
	E_1= \frac{\partial}{\partial t},
	\quad
	E_2=\frac{\partial}{\partial x},
	\quad
	E_3=t \frac{\partial}{\partial x}+ \frac{\partial}{\partial y},
	\quad
	E_4= \frac{t^2}{2} \frac{\partial}{\partial x} + t \frac{\partial}{\partial y} + \frac{\partial}{\partial z},
\end{equation*}
where $x,y,z$ are the coordinates on $\R^3$ and $t$ is the coordinate on $\R$ in $\R \ltimes \R^3$.

The group associated to $\mathrm{Nil}^4$ is $\mathrm{Nil}^4 \rtimes (\Z/2\Z)^2$ where $\mathrm{Nil}^4$ acts on itself by left translations and the copies of $\Z/2\Z$ are generated by the maps 
\begin{equation*}
  s_1: (t,x,y,z)\mapsto (t,-x,-y,-z), \quad 
  s_2: (t,x,y,z)\mapsto (-t,x,-y,z).
\end{equation*}

The homogeneous space $\mathrm{Nil}^4$ is also reductive, and thus we find a 4-parameter family of $G$-invariant Riemannian metrics,
\begin{equation*}
	\begin{aligned}
	&\tau_1 dt^2 + dx^2 + (t^2 + \tau_2) dy^2 + \left(\frac{t^4}{4} + t^2 (\alpha+\tau_2) + \tau_3 \right) dz^2 \\
	&- 2 t \,dx\,dy + (t^2 + 2\alpha) dx\,dz - (t^3 + 2t(\alpha + \tau_2) ) dy\,dz,
	\end{aligned}
\end{equation*}
where $\tau_1,\tau_2,\tau_3 \in \R^+$ and $\alpha \in \R\setminus\set{\pm\sqrt{\tau_3}}$.
Note that the frame
\[
	\left\{ \frac{E_1}{\sqrt{\tau_1}}, 
	E_2, 
	\frac{E_3}{\sqrt{\tau_2}}, 
	\frac{E_4-\alpha E_2}{\sqrt{\tau_3-\alpha^2}} \right\}
\]
is orthonormal with respect to these metrics.

In~\cite{wall-proceedings} it is shown that $\mathrm{Nil}^4$ does not admit any compatible complex structures.

\section{Future research on 4-dimensional Thurston geometries}
In this section we collect some ideas for future research in 4-dimensional Thurs\-ton geometries, mainly focussed on isometrically immersed submanifolds of these~spaces.
\begin{itemize}
	\item[(a)] Determine if the families of $G$-invariant metrics for $\mathrm{Nil}^4$ and $\mathrm{Sol}^4_1$ can be reduced, i.e.\ if it is possible to eliminate any of the parameters in these families.
	\item[(b)] In this article the Thurston geometry $\mathrm F^4$ was not described and no family of $G$-invariant metrics was exhibited.
	This is an obvious course for future research.
	\item[(c)] Study parallel and totally geodesic submanifolds in the Thurston geomtries $\mathrm{Sol}^4_1$, $\mathrm{Sol}^4_{m,n}$ and $\mathrm F^4$.
	Note that parallel and totally geodesic hypersurfaces of $\mathrm{Sol}^4_0$ are classified in~\cite{sol40} and of $\mathrm{Nil}^4$ in~\cite{nil4} for just one metric in the family of $G$-invariant metrics we exhibited in this article.
	\item[(d)] Formulate the fundamental equations for hypersurfaces in the 4-dimension\-al Thurston geometries with 4- or 5-dimensional isometry group.
	It might be an essential step to first develop a unified description of these spaces, in the same spirit as the  $\mathbb E(\kappa,\tau)$ spaces that Daniel used in~\cite{daniel2007isometric}.
\end{itemize}

\bibliographystyle{amsplain}
\bibliography{arxiv-version-ams}

\providecommand{\bysame}{\leavevmode\hbox to3em{\hrulefill}\thinspace}
\providecommand{\MR}{\relax\ifhmode\unskip\space\fi MR }
\providecommand{\MRhref}[2]{%
  \href{http://www.ams.org/mathscinet-getitem?mr=#1}{#2}
}
\providecommand{\href}[2]{#2}
\begin{thebibliography}{10}

\bibitem{wall-proceedings}
J.~C. Alexander, J.~L. Harer, and C.~T.~C. Wall, \emph{Geometries and geometric structures in real dimension 4 and complex dimension 2}, Geometry and Topology, Springer, 1985, pp.~268--292.

\bibitem{Backes-Reckziegel}
E.~Backes and H.~Reckziegel, \emph{On symmetric submanifolds of spaces of constant curvature}, Mathematische Annalen \textbf{263} (1983), 419--433.

\bibitem{Belkhelfa-Dillen-Inoguchi}
M.~Belkhelfa, F.~Dillen, and J.~Inoguchi, \emph{Surfaces with parallel second fundamental form in {B}ianchi-{C}artan-{V}ranceanu spaces}, Banach Center Publications \textbf{57} (2002), no.~1, 67--87.

\bibitem{bonahon}
F.~Bonahon, \emph{Geometric structures on 3-manifolds}, Handbook of geometric topology \textbf{93164} (2002).

\bibitem{daniel2007isometric}
B.~Daniel, \emph{Isometric immersions into 3-dimensional homogeneous manifolds}, Commentarii {M}athematici {H}elvetici \textbf{82} (2007), no.~1, 87--131.

\bibitem{sol40}
M.~D'haene, J.~Inoguchi, and J.~Van der Veken, \emph{Parallel and totally umbilical hypersurfaces of the four-dimensional {T}hurston geometry {$\mathrm{Sol}^4_0$}}, arXiv preprint 2303.14105 (2023).

\bibitem{nil4}
N.~Djellali, A.~Hasni, A.~M. Cherif, and M.~Belkhelfa, \emph{Classification of {C}odazzi and minimal hypersurfaces in {$\mathrm{Nil}^4$}}, arXiv preprint 2205.09848 (2022).

\bibitem{Filip}
R.~Filipkiewicz, \emph{Four-dimensional geometries}, Ph.D. thesis, University of Warwick, 1983.

\bibitem{hillman2002four}
J.~A. Hillman, \emph{Four-manifolds, geometries and knots}, vol.~5, University of Warwick, Mathematics Institute, 2002.

\bibitem{I-VdV-2007}
J.~Inoguchi and J.~Van~der Veken, \emph{Parallel surfaces in the motion groups {$E(1,1)$} and {$E(2)$}}, Bulletin of the {B}elgian {M}athematical {S}ociety-Simon Stevin \textbf{14} (2007), no.~2, 321--332.

\bibitem{I-VdV-2008}
\bysame, \emph{A complete classification of parallel surfaces in three-dimensional homogeneous spaces}, Geometriae Dedicata \textbf{131} (2008), 159--172.

\bibitem{Lawson}
H.~B. Lawson, \emph{Local rigidity theorems for minimal hypersurfaces}, Annals of Mathematics \textbf{89} (1969), no.~1, 187--197.

\bibitem{maier}
S.~Maier, \emph{Conformal flatness and self-duality of {T}hurston geometries}, Proceedings of the American Mathematical Society \textbf{126} (1998), no.~4, 1165--1172.

\bibitem{Takeuchi}
M.~Takeuchi, \emph{Parallel submanifolds of space forms}, Manifolds and Lie groups, Springer, 1981, pp.~429--447.

\bibitem{thurston-book}
W.~P. Thurston and S.~Levy, \emph{Three-dimensional geometry and topology}, Princeton mathematical series 35, Princeton university press, Princeton (N.J.), 1997 (eng).

\bibitem{Wall}
C.~T.~C. Wall, \emph{Geometric structures on compact complex analytic surfaces}, Topology \textbf{25} (1986), no.~2, 119--153.

\end{thebibliography}

\end{document}